\documentclass[12pt]{amsart}

\usepackage{latexsym}
\usepackage{amsmath}
\usepackage{amssymb}
\usepackage{multicol}

\newtheorem{theorem}{\bf Theorem}[section]

\newtheorem{proposition}[theorem]{\bf Proposition}
\newtheorem{lemma}[theorem]{\bf Lemma}

\newtheorem{definition-theorem}[theorem]{\bf Theorem-Definition}

\def\S{{\mathcal S}}

\def\i{{\rm i}}

\def\R{\mathbb{R}}
\def\C{\mathbb{C}}
\def\Z{\mathbb{Z}}

\def\g{\mathfrak{g}}
 \def\k{\mathfrak{k}}

\def\g{\frak{g}}

\def\p{\mathfrak{p}}
\def\a{\mathfrak{a}}

\def\k{\mathfrak{k}}

\include{include}

\begin{document}

\noindent {\LARGE \bf Equivariant cohomology of real flag manifolds}
\bigskip\\
{\bf Augustin-Liviu Mare } \\
{\small {\it Department of Mathematics and Statistics, University of Regina,
Regina SK, S4S 0A2 Canada\\
E-mail address:} {\tt mareal@math.uregina.ca}}
\smallskip \\
\title{}

\vspace{-1.5cm}

\begin{abstract}
Let $P=G/K$ be a semisimple non-compact Riemannian symmetric space,
where $G=I_0(P)$ and $K=G_p$ is the stabilizer of $p\in P$.  Let $X$
be an  orbit of the (isotropy) representation of $K$ on $T_p(P)$
($X$ is called a real flag manifold). Let $K_0\subset K$ be the
stabilizer of a maximal flat, totally geodesic submanifold of $P$
which contains $p$. We show that if all the simple root
multiplicities of $G/K$ are at least 2 then $K_0$ is connected and
the action of $K_0$ on $X$ is equivariantly formal. In the case when
the multiplicities are equal and at least 2, we will give a purely
geometric proof of a formula of Hsiang, Palais and Terng concerning
$H^*(X)$. In particular, this gives a conceptually new proof of
Borel's formula for the cohomology ring of an adjoint orbit  of a
compact Lie group.

\vspace{0.2cm}

 \noindent {\it
MS Classification}:  53C50, 53C35, 57T15

\end{abstract}

\maketitle




\section{Introduction}

Let  $G/K$ be a non-compact symmetric space, where $G$ is a
non-compact connected semisimple Lie group and $K\subset G$ a
maximal compact subgroup. Then $K$ is connected [He, Thm. 1.1, Ch.
VI] and there exists a Lie group automorphism $\tau$ of $G$ which is
involutive and whose fixed point set is $G^{\tau}=K$. The involutive
automorphism $d(\tau)_e$ of $\g={\rm Lie}(G)$ induces the Cartan
decomposition
$$\g =\k \oplus \p,$$ where $\k$ (the same as ${\rm Lie}(K)$) and $\p$ are the
$(+1)$-, respectively $(-1)$-eigenspaces of $(d\tau)_e$.   Since
$[\k,\p] \subset \p$, the space $\p$ is $Ad_G(K):=Ad(K)$-invariant.
The orbits of the action of $Ad(K)$ on $\p$ are called {\it  real
flag manifolds}, or {\it s-orbits}. The restriction of the Killing
form of $\g$ to $\p$ is an $Ad(K)$-invariant inner product on $\p$,
which we denote  by $\langle \ , \ \rangle$.

Fix $\a \subset \p$ a maximal abelian subspace. Recall that the
roots of the symmetric space $G/K$ are linear functions $\alpha
:\a\to \R$ with the property that the space
$$\g_{\alpha}:=\{z\in \g : [x,z]=\alpha(x)z \ {\rm for \ all }
\ x\in\a\}$$ is non-zero. The set $\Pi$ of all roots is a root
system in $(\a^*, \langle \ , \ \rangle)$. Pick $\Delta\subset
\Pi$ a simple root system and let $\Pi^+\subset \Pi$ be the
corresponding set of positive roots. For any $\alpha\in\Pi^+$ we
have
$$\g_{\alpha}+\g_{-\alpha}=\k_{\alpha}+\p_{\alpha},$$
where $\k_{\alpha}=(\g_{\alpha}+\g_{-\alpha})\cap\k$ and
$\p_{\alpha}=(\g_{\alpha}+\g_{-\alpha})\cap\p$. We have the direct
decompositions
$$\p=\a+\sum_{\alpha \in \Pi^+}\p_{\alpha},\quad \k=\k_0+\sum_{\alpha \in \Pi^+}\k_{\alpha},$$
where $\k_0$ denotes the centralizer of $\a$ in $\k$. The {\it
multiplicity}
 of a root $\alpha\in \Pi^+$ is
$$m_{\alpha}=\dim \k_{\alpha} +\dim \k_{2\alpha}.$$
We note that this definition is slightly different from
the standard one (see e.g. [Lo, Ch. VI, section 4]) which says that the multiplicity of
$\alpha$ is just $\dim\k_{\alpha}$.

Now $\k_0$ is the Lie algebra of the Lie group $K_0:=C_K(\a)$
 as well as of
$K_0':=N_K(\a)$. One can see that $K_0$ is a normal subgroup of
$K_0'$; the {\it Weyl group} of the symmetric space $G/K$ is
$$W=K_0'/K_0.$$ It can be realized geometrically as the (finite)
subgroup of $O(\a,\langle \ , \ \rangle)$ generated by the
reflections about the hyperplanes $\ker \alpha$, $\alpha\in\Pi^+$.

Take $x_0\in\a$ and let $X=Ad(K)x_0$ be the corresponding flag
manifold. The goal of our paper is to describe the cohomology, always with
coefficients in $\R$, of $X$. The first main result   concerns the action of
$K_0$ on $X$.

\begin{theorem}\label{firstmain}
If the symmetric space $G/K$ has all root multiplicities
$m_{\alpha}$, $\alpha\in\Delta$,
 strictly greater than 1 then:

(a) $K_0$ is connected;

(b) the action of $K_0$ on $X=Ad(K)x_0$ is equivariantly
formal, in the sense that $$H_{K_0}^*(X)\simeq H^*(X)\otimes
H_{K_0}^*({\rm pt})$$ by an isomorphism of $H_{K_0}^*({\rm
pt})$-modules;

(c) we have the  isomorphisms of $\R$-vector spaces
$$H^*(X)\simeq \sum_{w\in W}H^{*-d_w}(w.x_0), \quad
H^*_{K_0}(X)\simeq \sum_{w\in W}H_{K_0}^{*-d_w}(w.x_0).$$ Here
$$d_w=\sum m_{\alpha}$$ where the sum runs after all
$\alpha\in\Pi^+$ such that
 $\alpha/2\notin\Pi^+$  and   the
 line segment  $[x_0, w.x_0)$ crosses the hyperplane $\ker\alpha.$
\end{theorem}

{\bf Remark.}  Let $U$ be the (compact) Lie subgroup  of $G^{\C}$ whose Lie algebra is $\k \oplus
\i\p$. Then the manifold $X=Ad(K)x_0$ is the ``real locus" [Go-Ho], [Bi-Gu-Ho] of an anti-symplectic
involution on the adjoint orbit $Ad(U)\i x_0$ (see e.g. [Du, section 5]). 
The natural action of the torus $T:=\exp (\i\a)$  on this orbit is Hamiltonian. 
In this way, $X$ fits into
the more general framework  of [Go-Ho] and [Bi-Gu-Ho]. But these papers
investigate $X$ from the perspective of the action
of $T_{\R}=T\cap K=T\cap K_0$, whereas we
are interested here in the action on $X$ of a group which may be  larger than $T_{\R}$,  namely $K_0$.

In the second part of our paper we will deal with the ring structure of the usual
cohomology of $X$, under the supplementary assumption that   the
symmetric space has all root multiplicities equal. By [He, Ch. X, Table VI], their common value 
can be only 2, 4 or 8. An important
ingredient is the action of $W=K_0'/K_0$ on
$X$ given by
\begin{equation}\label{action}hK_0. Ad(k)x_0 = Ad(k) Ad(h^{-1})x_0,\end{equation} for any $h
\in K_0'$ and  $k\in K$. By functoriality, this induces an action of $W$ on $H^*(X)$. We also note that
$W$ acts in a natural way on $\a^*$.

\begin{theorem}\label{secondmain}
Assume that $G/K$ is an irreducible non-compact symmetric space whose simple root multiplicities are equal to
the same number, call it $m$, which is at least 2. Take $X=Ad(K)x_0$.

(i) If   $x_0$
is a regular point of $\a$, then there exists a canonical linear $W$-equivariant isomorphism
$\Phi: \a^*\to H^m(X)$. Its natural extension $\Phi:S(\a^*) \to H^*(X)$ is a surjective ring homomorphism
whose  kernel 
is the ideal $\langle S(\a^*)_+^W\rangle$ generated by all nonconstant $W$-invariant elements of
$S(\a^*)$. Consequently we have the $\R$-algebra isomorphism
$$H^*(X)\simeq S(\a^*)/\langle S(\a^*)_+^W\rangle.$$

(ii) If $x_0$ is an arbitrary point in $\a$, then we have the $\R$-algebra isomorphism
$$H^*(X)\simeq S(\a^*)^{W_{x_0}}/\langle S(\a^*)_+^W\rangle,$$
where $W_{x_0}$ is the $W$-stabilizer of $x_0$.

\end{theorem}



{\bf Remark.}  Any real flag manifold $X=Ad(K)x$ with the canonical
embedding in $(\p, \langle \ , \ \rangle )$ is an element of an
{\it isoparametric foliation} [Pa-Te]. The topology of such
manifolds, including their cohomology rings, has been investigated by Hsiang, Palais and Terng
in [Hs-Pa-Te] (see also [Ma]). 
The formulas for $H^*(X)$ given by Theorem \ref{secondmain} 
 have been proved by them in that paper. Even though we do use some of their ideas (originating in
[Bo-Sa]), our proof is different: they rely on Borel's  formula
[Bo] for the cohomology of a generic adjoint orbit of a compact Lie group,  whereas we actually prove it. 

{\bf Acknowledgements.} I  thank Jost Eschenburg for
discussions about the topics of the paper. I also thank
Tara Holm as well as    the 
referees for suggesting several improvements.

\section{Symmetric spaces with multiplicities at least 2 and
their $s$-orbits}

Let $G/K$ be an arbitrary non-compact symmetric space, $x_0\in\a$
and $X=Ad(K)x_0$ the corresponding $s$-orbit. The latter is a
submanifold of the Euclidean space $(\p,\langle \ ,\ \rangle)$.
The Morse theory of height functions on $X$ will be an essential
instrument. The following proposition summarizes results from
[Bo-Sa] or [Hs-Pa-Te] (see also [Ma]).

\begin{proposition}
(i) If $a\in \a$ is a general vector (i.e. not contained
in any of the hyperplanes $\ker\alpha$, $\alpha\in\Pi^+$), then
the height function $h_a(x)=\langle a,x\rangle$, $x\in X$ is a
Morse function. Its critical set is the orbit $W.x_0$.

(ii) Assume that $a$ and $x_0$ are contained in the same
Weyl chamber in $\a$. Then the index of $h_a$ at the critical
point $w.x_0$ is \begin{equation}\label{index}d_w=\sum
m_{\alpha}\end{equation} where the sum runs after all
$\alpha\in\Pi^+$ such that
 $\alpha/2\notin\Pi^+$  and   the
 line segment  $[a, wx_0)$ crosses the hyperplane $\ker\alpha.$
\end{proposition}

In the  next lemma we consider the situation when all 
root multiplicities are at least 2.

\begin{lemma}\label{cor1} Assume that the root multiplicities $m_{\alpha}$,
$\alpha\in \Delta$, of the symmetric space $G/K$ are all strictly
greater than 1. Then:

(i) for any general vector $a\in\a$, the height function
$h_a:X\to \R$ is $\Z$-perfect,

(ii) the space $K_0$ is connected,

(iii) if $X=Ad(K)x_0$, then the orbit $W.x_0$ is
contained  in the fixed point set $X^{K_0}.$
\end{lemma}

\begin{proof}
(i) According to [Ko, Theorem 1.1.4], there exists a metric on $X$ such that if two critical points
$x$ and $y$  
can be joined by a gradient line, then $x=s_{\gamma}y$, where $\gamma\in \Pi^+$.
 By (\ref{index}), the difference of the indices of
$x$ and $y$  is different from $\pm 1$. 
Because the stable and unstable manifolds intersect transversally  [Ko, Corollary 2.2.7], the 
Morse complex of $h_a$ has all boundary operators identically zero, hence $h_a$ is $\Z$-perfect.

 (ii) Take $a\in\a$ a
general vector. The height function $h_a$  on $Ad(K)a$  is
$\Z$-perfect. From  (\ref{index}) we deduce that $H_1(Ad(K)a,
\Z)=0$, thus $Ad(K)a$ is simply connected. On the other hand, the
stabilizer $C_K(a)$ is just $K_0$ (see e.g. [Bo-Sa]). Because  $K/K_0$ is
simply connected and $K$ is connected, we deduce that $K_0$ is
connected.

(iii) The height function $h_a$ is $Ad(K_0)$-invariant, thus
$Crit(h_a)=W.x_0$ is also $Ad(K_0)$-invariant. The result follows
from the fact that $K_0$ is connected. \hfill{{$\square$}}

\end{proof}

We are now ready to prove Theorem \ref{firstmain}.

\noindent {\it Proof of Theorem \ref{firstmain}}. Point (a) was proved in Lemma \ref{cor1} (ii).

(b) According to
[Gu-Gi-Ka, Proposition C.25] it is sufficient to show that
$H^*_{K_0}(X)$ is free as a $H^*_{K_0}({\rm pt})$-module. In order to do
that we  consider the height function $h_a:X\to\R$ corresponding
to a general $a\in\a$. We use the same arguments as in the proof of Lemma \ref{cor1}, (i).
The  function $h_a$  is   a $K_0$-invariant. By the same reasons as above, 
the $K_0$-equivariant Morse
complex [Au-Br, Sections 5 and 6] has all boundary operators identically
zero.  Thus
$H^*_{K_0}(X)$ is a free $H^*_{K_0}({\rm pt})$-module (with a basis
indexed by ${\rm Crit} (h_a)=W.x_0$).


(c)  The space $H_*(X)$ has a basis $\{[X_{w.x_0}]:w\in W\}$,
where ${X}_{w.x_0}$ is some $d_w$-dimensional cycle in $X$, $w\in
W$. The evaluation pairing $H^*(X)\times H_*(X)\to \R$ is
non-degenerate; consider the basis of $H^*(X)$ dual to
$\{[{X}_{w.x_0}]:w\in W\}$, which gives one element of degree
$d_w$ for each $w.x_0$. The result follows.\hfill{{$\square$}}



\section{Cohomology of $s$-orbits of symmetric spaces with uniform
multiplicities at least 2}

Throughout this section  $G/K$ is a non-compact irreducible
symmetric space whose simple root multiplicities are all equal to
$m$,   where $m\ge 2$;  $x_0\in\a$ is a regular element and
$$X=Ad(K)x_0\simeq K/K_0$$ is the corresponding real flag manifold.
There are three such symmetric spaces; their compact duals are (see e.g. [Hs-Pa-Te, Section 3]):
\begin{itemize}
\item[1.] any connected simple compact Lie group $K$; we have
$m=2$; the flag manifold is $X=K/T$, where $T$ is a maximal torus
in $K$;
\item[2.] $SU(2n)/Sp(n)$ where $m=4$; the flag manifold is
$X=Sp(n)/Sp(1)^{\times n}$;
\item[3.] $E_6/F_4$ where $m=8$; the flag manifold is
$X=F_4/Spin(8)$.
\end{itemize}

Let $\Delta=\{\gamma_1, \ldots, \gamma_l\}$ be a simple root
system of $\Pi$. To each $\gamma_j$ corresponds the distribution
$E_j$ on $X$, defined as follows: its value at $x_0$ is
$$E_j  (x_0) =[\k_{\gamma_j},x_0]$$
and  $E_j$ is $K$-invariant, i.e.
$$E_j(Ad(k)x_0)=Ad(k)E_j(x_0),$$
for all $k\in K$.

A basis of $H_m(X)$ can be obtained as follows: Assume that $x_0$ is in
the (interior of the) Weyl chamber $C\subset \a$ which is bounded
by the hyperplanes $\ker\gamma_j$, $1\le j\le l$. The Weyl group
$W$ is generated by $s_j$, which is the reflection of $\a$ about the wall
$\ker\gamma_j$, $1\le j\le l$. For each $1\le j \le l$ we consider
the Lie subalgebra $\k_0+\k_{\gamma_j}$ of $\k$; denote by $K_j$ the
corresponding connected subgroup  of $K$. It turns out that the
orbit $Ad(K_j)x_0$ is a round $m$-dimensional metric sphere in
$(\p, \langle \ , \ \rangle)$. To any $x=Ad(k)x_0\in X$ we attach
the round sphere
$$S_j(x)=Ad(k)Ad(K_j)x_0.$$
The spheres $S_j$ are integral manifolds of the distribution
$E_j$. We denote by $[S_j]$ the homology class carried by any of
the spheres $S_j(x)$, $x\in X$. It turns out that $S_1(x_0),
\ldots, S_l(x_0)$ are cycles of Bott-Samelson type (see [Bo-Sa], [Hs-Pa-Te])
 for the index $m$ critical points of the height function
$h_a$, thus $[S_1],\ldots, [S_l]$ is a basis of $H_m(X)$.

 The following result concerning the action of $W$
on $H_m(X)$ was proved in [Hs-Pa-Te, Corollary 6.10] (see also
[Ma, Theorem 2.1.1]):

\begin{proposition}\label{homology} We can  choose an orientation
of the spheres $S_j$, $1\le j\le l$, such that the linear
isomorphism $\a\to H_m(X)$ determined by
$$\gamma_j^{\vee}:=
\frac{2\gamma_j}{\langle \gamma_j,\gamma_j\rangle}\mapsto [S_j],$$
$1\le j\le l$, is $W$-equivariant.
\end{proposition}

We need one more result concerning the action of $W$ on $H^*(X)$:

\begin{lemma}\label{zero}
Let $x\in\a$ be an arbitrary element, $C=C_K(x)$  its centralizer in $K$,  and 
let $$p:X=K/K_0 \to Ad(K)x=K/C$$ be the natural
map induced by the inclusion $K_0\subset C$. Then the map $p^*:H^*(Ad(K)x)\to H^*(X)$ is
injective. Its image is
$$p^*H^*(Ad(K)x) = H^*(X)^{W_{x}}$$
where the right hand side denotes the set of all $W_x$-invariant
elements of $H^*(X)$. Here $W_x$ denotes the $W$-stabilizer of
$x$. In particular, the only elements in $H^*(X)$ which are
$W$-invariant are those of degree $0$, i.e.
$$H^*(X)^W=H^0(X).$$

\end{lemma}

\begin{proof} The map $p:K/K_0\to K/C$ is  a fibre bundle.
The fiber  $C/K_0$ is an $s$-orbit of the symmetric space $C_G(x)/C_K(x)$.
The latter has  all root multiplicities equal to $m$, as they are
all  root multiplicities of some roots of $G/K$. By Theorem
\ref{firstmain} (ii), $C/K_0$ can have non-vanishing cohomology groups
only in dimensions which are multiples of $m$. The same can be
said about the cohomology of the space $K/C$. Because
$m\in \{ 2,4,8\}$,  the spectral sequence of the bundle $p:K/K_0\to
K/C$ collapses, which implies that $p^*$ is injective.

The map $p$ is $W$-equivariant with respect to the actions of $W$ on
$Ad(K)x_0$, respectively $Ad(K)x$ defined by (\ref{action}). Thus
if $w\in W_x$, then $w|_{Ad(K)x}$ is the identity map, hence we
have $p\circ w=p$. This implies the inclusion
$$p^*H^*(Ad(K)x)\subset H^*(X)^{W_x}.$$

On the other hand, the action of $W$ on $X$ defined by
(\ref{action}) is free, as the $Ad(K)$ stabilizer of the general
point $x_0$ reduces to $K_0$. Consequently we have
$$H^*(X)^{W_x}=H^*(X/W_x)$$
and $$\chi(X/W_x)=\frac{\chi(X)}{|W_x|}=\frac{|W|}{|W_x|},$$
where $\chi$ denotes the Euler-Poincar\'e characteristic. It
follows from Theorem \ref{firstmain} (c) that
$$\dim H^*(X)^{W_x} = \frac{|W|}{|W_x|} =\dim H^*(Ad(K)x).$$
Now we use that $p^*$ is injective.

In order to prove the last statement of the lemma, we take
$x=0\in\a$.\hfill{{$\square$}}

\end{proof}


Let us consider  the Euler class $\tau_i=e(E_i)\in H^m(X)$, $1\le
i\le l$. We will prove that:
\begin{lemma}\label{euler}
(i) The cohomology classes $\tau_i$, $1\le i\le l$ are a basis of
$H^m(X)$.

(ii)   The linear isomorphism $\Phi: \a^*\to H^m(X)$ determined by
$$\gamma_i\mapsto e(E_i),$$
$1\le i\le l$, is $W$-equivariant.

\end{lemma}
\begin{proof} By Proposition \ref{homology} we know that
$$s_{i*}[S_j] =[S_j]-d_{ji}[S_i],$$
where $$d_{ji}= 2\frac{\langle
\gamma_j^{\vee},\gamma_i^{\vee}\rangle}{\langle
\gamma_i^{\vee},\gamma_i^{\vee}\rangle}.$$ Denote by $\langle \ , \ \rangle$ the evaluation pairing 
$H^2(M)\times H_2(M) \to \R$. Consider $\alpha_j\in H^2(M)$ such that $\langle \alpha_j, [S_i]\rangle =\delta_{ij}$,
$1\le i,j\le l$. 
 Take the expansion
$$\tau_i=\sum_{j=1}^lt_{ij}\alpha_{j}.$$
The automorphism $s_i$ of $X$ maps the distribution $E_i$ onto
itself and changes its orientation (since so does the antipodal
map on an $m$-dimensional sphere). Thus
$$s_i^*(\tau_i)=-\tau_i.$$
Consequently we have
$$t_{ij}=\langle \tau_i, [S_j] \rangle =
\langle -s_i^*(\tau_i), [S_j]\rangle = -\langle \tau_i,
s_{i*}[S_j]\rangle =-\langle \tau_i, [S_j]-d_{ji}[S_i]\rangle =
-t_{ij}+2d_{ji}$$ which implies $t_{ij}=d_{ji}$. By Proposition
\ref{homology}, the matrix $(d_{ij})$ is the Cartan matrix of the
root system dual to $\Pi$, hence it is non-singular. Consequently
$\tau_i$, $1\le i\le l$ is a basis of $H^m(X)$. Again by
Proposition \ref{homology} we have
$$\langle s_{j}^*(\tau_i) ,[S_k]\rangle=
\langle\tau_i, [S_k]-d_{kj}[S_j]\rangle=t_{ik}-d_{kj}t_{ij}
=t_{ik}-t_{jk}d_{ji}, $$ thus $$s_j^*(\tau_i) = \tau_i
-d_{ji}\tau_j.$$ It remains to notice that $d_{ji}$ can also be
expressed as
   $$d_{ji}=2\frac{\langle \gamma_i,\gamma_j\rangle}{\langle
\gamma_j,\gamma_j\rangle}.$$ \hfill{{$\square$}}\end{proof}

We are now ready to prove
Theorem \ref{secondmain}:

{\it Proof of Theorem \ref{secondmain}} 
(i) Consider the ring homomorphism $\Phi: S(\a^*)\to H^*(X)$ induced by
$\gamma_i\mapsto e(E_i)$, $1\le i\le l$.
By Lemma \ref{euler}, $\Phi$ is $W$-equivariant
and from Lemma \ref{zero} we deduce that $\langle S(\a^*)_+^W\rangle \subset \ker \Phi$. By  Lemma
\ref{last} (see below), it is sufficient to prove that
$$\Phi(\prod_{\alpha\in \Pi^+}\alpha)\ne 0.$$

To this end, we will describe explicitly $\Phi(\alpha)$, for
$\alpha\in \Pi^+$. Write $\alpha= w.\gamma_j$, where $w\in W$. The
latter is of the form $w=hK_0$, with $h\in K_0'$. The image of
$S_j(x_0)$ by the automorphism $w$ of $X$ is
\begin{align*}
w(S_j(x_0))=& Ad(K_j)Ad(h^{-1})x_0=Ad(h^{-1})Ad(hK_jh^{-1})x_0\\=&
Ad(h^{-1})Ad(K_{\alpha})x_0=Ad(h^{-1})S_{\alpha}(x_0)=
S_{\alpha}(Ad(h^{-1})x_0)\\=&S_{\alpha}(w.x_0).\end{align*} Here
$K_{\alpha}$ is the connected subgroup of $K$ of Lie algebra
$\k_0+\k_{\alpha}$ and $S_{\alpha}(x_0):= Ad(K_{\alpha})x_0$ is a
round metric sphere through $x_0$; for any $x=Ad(k)x_0\in X$ we
have  $S_{\alpha}(x):=Ad(k)S_{\alpha}(x_0)$, which is an integral
manifold of
$$E_{\alpha}(x)=Ad(k)[\k_{\alpha}, x_0].$$
It is worth mentioning  in passing that the spheres $S_{\alpha}$ and the distributions
$E_{\alpha}$ are the curvature spheres, respectively 
curvature distributions of the isoparametric submanifold
$X\subset \p$ (see the remark following Theorem \ref{secondmain} in the introduction).  Thus the differential
of $w$ satisfies $(dw)(E_j) =E_{\alpha}$, which implies
$$e(E_j) =w^* e(E_{\alpha}).$$
 Consequently
$$\Phi(\alpha) = \Phi(w.\gamma_j)=w^{-1}.\Phi(\gamma_j)=
(w^{-1})^*(e(E_j)) =e(E_{\alpha}).$$

We deduce that $$\Phi(\prod_{\alpha\in \Pi^+}\alpha)=
\prod_{\alpha\in \Pi^+} e(E_{\alpha}) = e(\sum_{\alpha\in \Pi^+}
E_{\alpha}). $$ On the other hand,
$$\sum_{\alpha\in \Pi^+} E_{\alpha}(x_0)=
\sum_{\alpha\in \Pi^+} [\k_{\alpha}, x_0]=[\k, x_0]= T_{x_0}X$$
thus $$\sum_{\alpha\in \Pi^+} E_{\alpha}= TX.$$ It follows that
$$\Phi(\prod_{\alpha\in \Pi^+}\alpha)=e(TX),$$
which is different from zero, as
$$e(TX)([X])=\chi(X)=|W|,$$
where $\chi(X)$ is the Euler-Poincar\'e characteristic of
$X$.

(ii) We apply Lemma \ref{zero}. 
\hfill{{$\square$}}

The following lemma has been used in the proof:

\begin{lemma}\label{last}{\rm ([Hi, Lemma 2.8])} Let $I$ be a graded ideal of $S(\a^*)$ which
is also a vector subspace and such that $\langle S(\a^*)_+^W\rangle\subset I$. We have $I=\langle S(\a^*)_+^W\rangle$
if and only if $$\prod_{\alpha\in \Pi^+}\alpha\notin I.$$

\end{lemma}

A proof of this lemma can also be found in the appendix. 



{\small

\section{Appendix: Proof of Lemma \ref{last}}

The goal of this appendix is to provide a proof of Lemma \ref{last}, which is stated
without a proof in  [Hi]. As mentioned in the introduction,
the Weyl group $W$ can be realized as the group of orthogonal transformations of $\a$ generated
by the reflections $s_{\alpha}$, $\alpha \in \Pi^+$. In fact, if $\{\gamma_1,\ldots,\gamma_l\}$
is a simple root system, then $W$ is generated by $s_i:=s_{\gamma_i}$, $1\le i\le l$.
Denote by $w_0$ the longest element of $W$, where the length is measured with respect to the
generating set $\{s_1,\ldots, s_l\}$.  We will use the notations
$$\S := S(\a^*), \quad I_W:= \langle S(\a^*)^W_+\rangle.$$
First of all we note that the action of $W$   on the polynomial ring $\S$ is given  by $$(w.f)(x)=f(w^{-1}.x),$$ 
where $w \in W$, $f
\in \S$, $x \in \a$. This action preserves the
grading of $\S$, hence the ideal $I_W$ generated by the
nonconstant $W$-invariant polynomials is also graded. The most prominent example of a
polynomial which is not $W$-invariant is
$$d=\prod _{\alpha \in \Pi ^+}\alpha.$$  In fact $d$ is
skew-invariant, in the sense that $w.d =(-1)^{l(w)}d$, for any $w
\in W$. 

If $\alpha\in \Pi^+$, we consider the operator $\Delta _{\alpha} :
\S \rightarrow \S$ defined as follows:
$$\Delta _{\alpha} (f) =\frac{f -s_{\alpha}.f}{\alpha},$$
$f \in \S$. Note that $f -s_{\alpha}.f$ vanishes on the space $\ker \alpha $,
hence $\Delta_{\alpha}(f)$ is really a polynomial.
The following result is straightforward:

\begin{lemma} If $w \in W$, $\alpha \in \Pi^+$, $f,g \in \S$, then we have:

(a) $\Delta _{\alpha}(fg)=\Delta _{\alpha} (f) g +s_{\alpha}(f)
\Delta _{\alpha}(g)$;

(b) $\Delta _{\alpha}(I_W) \subset I_W$.

\end{lemma}
To any $w\in W$ we can associate  the operator
$\Delta _{w} :\S \rightarrow \S$, which has degree $-l(w)$, and is
defined as follows:
take $w =s_{i_1} \ldots s_{i_k}$  a reduced expression and put
$\Delta _w= \Delta _{\gamma _{i_1}} \cdots \Delta_{\gamma _{i_k}}$.
We note that $\Delta_w$ does not depend on the choice of the reduced expression 
(see e.g. [Hi, Proposition 2.6]). 
The operators obtained in this way have the following property (see [Hi, Lemma 3.1]):
\begin{equation}\label{delta} \Delta _w  \circ \Delta _{w'} = \begin{cases}
                                        \Delta_{ww'},\ {\rm if} \ 
                                        l(ww')=l(w)+l(w')\\
                                        0,\ {\rm otherwise}
                                        \end{cases}
\end{equation}


A classical result which goes back to Chevalley, says that the ideal
$I_W$ is generated by $l$ homogeneous polynomials, which are algebraically
independent. Let $d_1,\cdots, d_l$ denote their degrees.
It follows that the Poincar\'e polynomial of $\S/I_W$ is:
$$P(\S/I_W)= \sum_{k=0}^{\infty}(\dim\S^k-
\dim I_W ^k)t^k= \prod_{j=1}^l (1+t+\cdots +t^{d_j-1}).$$
Combined with the fact that $d_1+\cdots d_l=N+l$
(see for instance [Hu, Theorem 3.9]), this tells us that
$I^k=\S^k$, for $k \geq N+1$. The same polynomial can be expressed as (see [Hu, Theorem 3.15]):
$$P(S/I_W)= \sum _{w\in W}t^{l(w)}.$$
We deduce that $\dim\S^k -\dim
I_W^k$ equals the number of $w \in W$ with $l(w)=k$, $0
\leq k \leq N$. The following result describes a direct complement of $I_W^k$ in $\S^k$:

\begin{proposition}\label{elements}
For any $0\leq k \leq N$, the elements $\Delta _w(d)$,
$w \in W$, $l(w)=N-k$ are linearly independent and span a direct
complement of $I_W^k$ in ${\mathcal S}^k$.
\end{proposition}

\begin{proof}
The number of elements of $W$ of length $k$ equals the number of elements
of length $N-k$, hence we only have to prove that the polynomials
$\Delta _w(d)$, where $l(w)=N-k$ are linearly independent and  their span
intersected with  $I_W$ is $\{0\}$.
To this end, it is sufficient to show that if $$\sum_{l(w)=N-k}
\lambda _w \Delta _w (d) \in I_W^k$$ then all $\lambda _w$ must vanish.
Indeed, if we fix $v\in W$ with $l(v)=N-k$, then by (\ref{delta}), we have 
$$\Delta_{w_0 v^{-1}}(\sum_{l(w)=N-k}
\lambda _w \Delta _w (d) )=\lambda_v .$$ The left hand side of this equation 
is in $I_W^0$, hence it must be 0.  

\end{proof}

We are ready to prove Lemma \ref{last}:

{\it Proof of Lemma \ref{last}} We prove by induction on $k$ that $I_W^k =I^k$,
$0 \leq k \leq N$. Things are clear for $k=N$: $I_W^N$ equals $I^N$ because 
$I_W^N\subset I^N\neq \S^N$ and the codimension of $I_W^N$ in $\S^N$ is 1
(see   Proposition \ref{elements}).
Now, from $I^{k+1} =I_W^{k+1}$ we deduce that $ I^{k} =I_W^{k}$. Suppose
that we have  $$f
:=\sum_{l(w)=N-k} \lambda _w \Delta _w (d)\in I^k,$$ where $\lambda _w \in \R$, not
all of them equal to 0.
 We will prove by induction on $m \in \{0,\ldots ,k\}$ the following claim

{\it Claim.} For any
$h_m \in \S^m$ and any  $\alpha _1, \ldots , \alpha _m \in
\Pi^+$, we have $$h_m \Delta _{\alpha _1}\circ \ldots \circ
\Delta _{\alpha _m}(f) \in I^k.$$ 

For $m=0$, this is trivial.
Suppose it is true for a certain  $m$ and prove it for $m+1$. If $h_m
\in \S_m$, $\alpha _1, \ldots ,\alpha _m \in \Pi ^+$, $h$
an arbitrary homogeneous polynomial of degree 1, and $\alpha$ a
 positive root, then we have $$hh_m \Delta _{\alpha_1}\circ
\ldots \circ \Delta _{\alpha _m} (f) \in I^{k+1}=I_W^{k+1},$$ hence
its image by $\Delta _{\alpha}$ is in $I_W^k \subseteq I^k$.
We deduce that $$\Delta _{\alpha} (h) h_m \Delta _{\alpha_1}\circ
\ldots \circ  \Delta _{\alpha_m}(f)+s_{\alpha}(h)\Delta
_{\alpha}(h_m) \Delta _{\alpha_1} \circ \ldots \circ \Delta
_{\alpha_m}(f)+ s_{\alpha}(h)s_{\alpha}(h_m)\Delta _{\alpha} \circ
\Delta _{\alpha_1} \circ \ldots \circ\Delta _{\alpha_m} (f)$$
is in $I^k$, consequently $s_{\alpha}(hh_m)  \Delta
_{\alpha}\circ \Delta _{\alpha_1} \circ \ldots \circ \Delta
_{\alpha_m} (f) \in I^k$. Since any $h_{m+1} \in \S^{m+1}$
is a linear combination of polynomials of the form $s_{\alpha}(hh_m)
$, the claim is proved.

We deduce that for any $v \in W$ with $l(v) =k$, and any
$h_k \in \S^k$ we have that $$h_k \Delta _{v}(f) \in I^k.$$
Fix now $w \in W$ with $l(w)=N-k$ and take $v:= w_0 w^{-1}$. Then
$\Delta_{v}(f) =\lambda _w $ by (1), hence $\lambda _w h_k
\in I^k$, for any $h_k \in \S^k$. But then $\lambda _w$
must vanish, since $I^k\neq \S^k$ (if they were equal, from $k\le N$ we would deduce
$I^N=S^N$, which is false). We conclude that $f=0$, which is a contradiction.
This finishes the proof. \hfill $\square$

}

\newpage

\bibliographystyle{abbrv}

\end{document}